\documentclass[12pt,dvips]{article}
\usepackage{rotating}
\usepackage{axodraw}
\usepackage{epsfig}
\usepackage{color}
\usepackage{cite}

\begin{document}

\def\thefootnote{\fnsymbol{footnote}}

\begin{figure*}[t]
\begin{flushleft}
\resizebox{2cm}{!}{\includegraphics{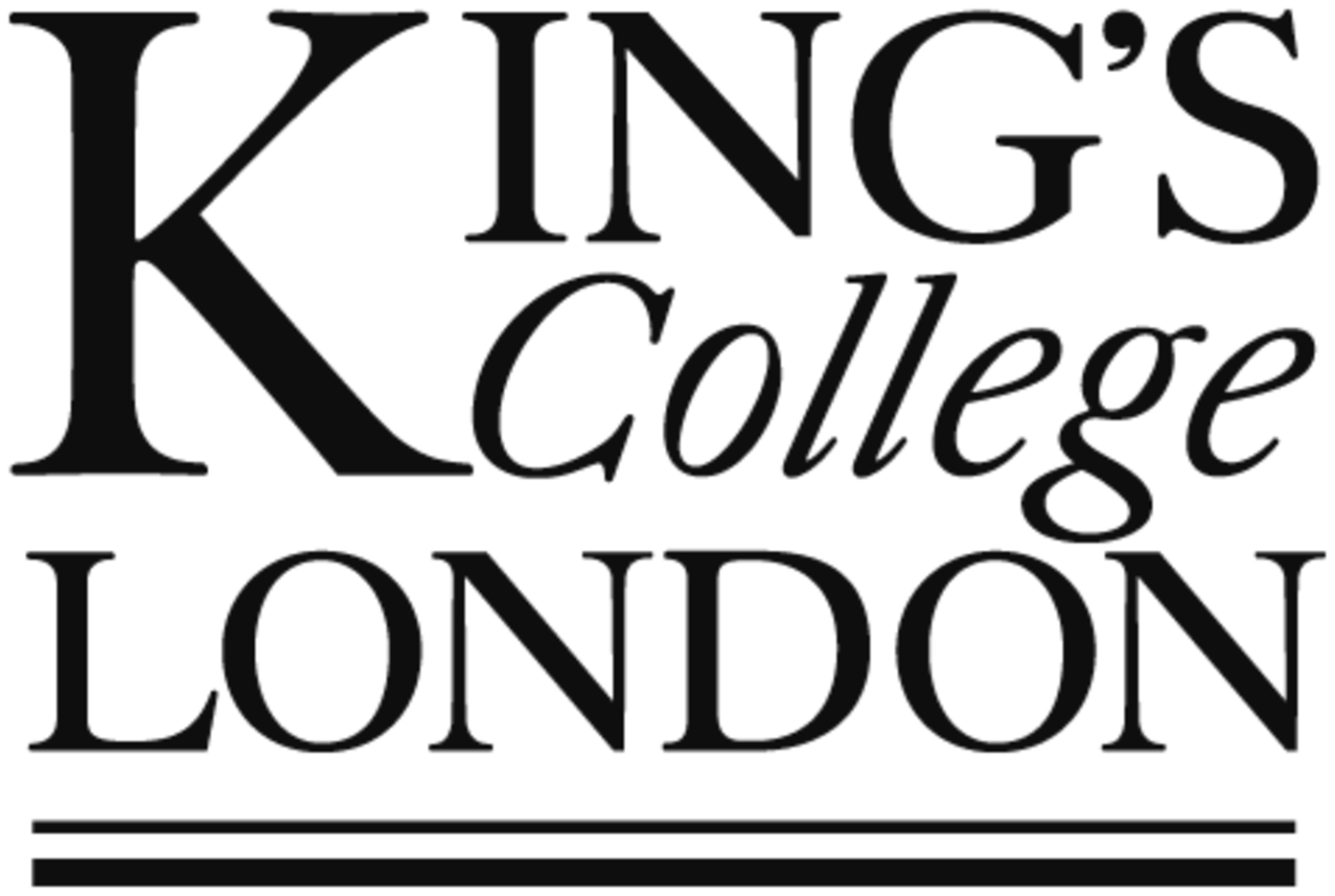}}
\end{flushleft}
\end{figure*}
${}$~\vspace{-4.8cm}
\begin{flushright}
{\tt CERN-PH-TH/2010-199}, {\tt KCL-PH-TH/2010-19}, {\tt MAN/HEP/2010/18 }\\
{~}\\
September 2010
\end{flushright}

\vspace{1cm}
\begin{center}
{\bf {\Large Note on a Differential-Geometrical Construction of \\
\vskip 0.1in
Optimal Directions in Linearly-Constrained Systems}}
\end{center}

\medskip

\begin{center}{\large
John~Ellis$^{a,b}$,
Jae~Sik~Lee$^c$ and
Apostolos~Pilaftsis$^d$}
\end{center}

\begin{center}
{\em $^a$Theory Division, CERN, CH-1211 Geneva 23,
  Switzerland}\\[0.2cm]
{\em $^b$Theoretical Physics and Cosmology Group, Department of
  Physics, King's~College~London, London WC2R 2LS, United Kingdom}\\[0.2cm] 
{\em $^c$Physics Division, National Center for Theoretical Sciences, 
Hsinchu, Taiwan 300}\\[0.2cm]
{\em $^d$School of Physics and Astronomy, University of Manchester,}\\
{\em Manchester M13 9PL, United Kingdom}
\end{center}

\bigskip\bigskip

\centerline{\bf ABSTRACT}

\noindent  
This note  presents an analytic construction of  the optimal unit-norm
direction ${\bf \widehat{x}} =  {\bf x}/||{\bf x}||$ that maximizes or
minimizes  the   objective  linear  expression,   ${\bf  B}\cdot  {\bf
\widehat{x}}$, subject to  a system of linear constraints  of the form
$[{\bf A}]\cdot  {\bf x}  = {\bf  0}$, where ${\bf  x}$ is  an unknown
$n$-dimensional  real  vector  to  be  determined  up  to  an  overall
normalization constant,  ${\bf 0}$ is an  $m$-dimensional null vector,
and  the  $n$-dimensional  real  vector  ${\bf B}$  and  the  $m\times
n$-dimensional real matrix~$[{\bf A}]$  (with $0\le m < n$) are given.
The analytic solution  to this problem can be expressed  in terms of a
combination of  double wedge  and Hodge-star products  of differential
forms.

\medskip
\noindent
{\small {\sc Keywords}: Linear programming; Linear algebra;
  Differential forms\\[3mm]
{\sc AMS Classification 2010}: 90C05, 15A69, 58A10.} 

\newpage

\section{Introduction}

A basic  problem in linear programming~\cite{JKS}  is the optimization
of the objective function, ${\bf  B}\cdot {\bf x}$, which is linear in
an   $n$-dimensional  real   vector  ${\bf   x}$,   with  non-negative
components, that is  subject to a system of  linear constraints of the
form $[{\bf  A}]\cdot {\bf x}  = {\bf A_0}$,  where $[{\bf A}]$  is an
$m\times n$-dimensional real matrix,  with $0\le m<n$, and ${\bf A_0}$
is an $m$-dimensional real  vector with non-negative components.  

This  note presents  an  exact,  analytic solution  to  a simple,  but
different class  of linear programming problems, where  ${\bf A}_0$ is
the $m$-dimensional null vector ${\bf  0}$ and the objective is ${\bf
B}\cdot {\bf \widehat{x}}$, where  ${\bf \widehat{x}} = {\bf x}/||{\bf
x}||$ is  the optimal  unit-norm direction to  be determined.   In the
absence of  any linear constraint  ($m=0$), the solution  is trivially
given  by  ${\bf  \widehat{x}}  =  {\bf  B}/||{\bf  B}||  \equiv  {\bf
\widehat{B}}$.  For a large number $m$ of linear constraints, however,
the   solution  becomes   non-trivial  and   relies  on   a  geometric
construction expressed  in terms of  wedge and Hodge-dual  products of
differential forms~\cite{Hestenes,SG,SMC}.

As an  example of this type  of linear programming  problem, this note
first describes the simple case involving a single constraint ${\bf A}
\cdot  {\bf x}  = 0$  on a  3-dimensional real  vector~${\bf  x}$. The
subspace  of  parameters  satisfying   such  a  linear  constraint  is
characterized, and  the analytic construction  of a real  vector ${\bf
  x}$  that  optimizes  the  linear  expression,  ${\bf  B}\cdot  {\bf
  \widehat{x}}$, is then described.  Subsequently, a generalization to
multiple constraints in higher-dimensional spaces is described.
We emphasize the formal analogies between the three- and 
higher-dimensional cases.

\section{A Three-Dimensional Example}

Let us  first consider a simple three-dimensional  (3D) example, where
${\bf x}  = (  x_1,\, x_2,\,  x_3 )$. Using  the standard  notation of
inner product multiplication  between vectors, the aim is  to find the
unit-norm direction  ${\bf \widehat{x}} = {\bf  x}/||{\bf x}||$, along
which the  expression, ${\bf B}\cdot {\bf  \widehat{x}}$, is maximized
or minimized,  subject into the constraint:  ${\bf A} \cdot  {\bf x} =
0$, where the three-vectors ${\bf A}$ and ${\bf B}$ are given.

\begin{figure}[t]
\begin{center}
\begin{picture}(200,180)(0,50)
\SetWidth{0.9}

\LongArrow(100,100)(200,80)\Text(205,80)[l]{\boldmath$\widehat{x}_2$}
\LongArrow(100,100)(100,185)\Text(100,190)[b]{\boldmath $\widehat{x}_3$}
\LongArrow(100,100)(30,60)\Text(25,60)[r]{\boldmath $\widehat{x}_1$}
\Line(100,100)(170,140)

\LongArrow(100,100)(130,150)\Text(125,155)[lb]{${\bf A}$}
\LongArrow(100,100)(70,150)\Text(80,155)[rb]{${\bf B}$}
\LongArrow(100,100)(70,130)\Text(70,125)[rt]{${\bf x}(t)$}

\Line(30,110)(100,150)
\DashLine(85,50)(150,87){3}
\Line(30,110)(57,80)\DashLine(57,80)(85,50){3}
\Line(100,150)(127,120)\DashLine(127,120)(155,90){3}

\end{picture} 
\end{center}
\caption{\it Three dimensional  example that illustrates the geometric
  construction of  the optimal direction  for ${\bf B}\cdot  {\bf x}$,
  subject to the  constraint ${\bf A}\cdot {\bf x}  = 0$.  The optimal
  direction  is  given  by  the  vector ${\bf  x}(t)$,  which  is  the
  intersection  of the  indicated  plane perpendicular  to the  vector
  ${\bf A}$ with the plane defined  by the vectors ${\bf B}$ and ${\bf
  A}$.}
\label{fig:geometry}
\vspace{1cm}
\end{figure}
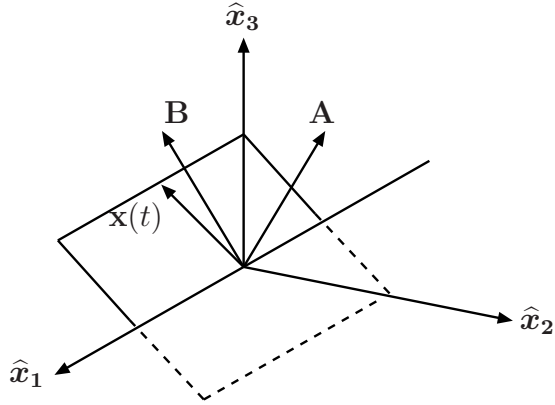

As represented schematically  in Fig.~\ref{fig:geometry}, the solution
to this  problem is  well known  and is given  by the  ordinary vector
triple product:
\begin{equation}
  \label{Param3}
{\bf x}(t)\ =\ t\ {\bf A}\times \Big( {\bf B} \times {\bf A} \Big)
\ ,
\end{equation}
where $t \in \mbox{I$\!$R}$ is an arbitrary real parameter. Evidently,
the optimal (maximal) direction is given by the unit-norm vector,
\begin{equation}
  \label{Optimal3}
{\bf \widehat{x}}\ =\ \frac{{\bf A}\times \Big( {\bf B} \times
  {\bf A} \Big)}{||{\bf A}\times \Big( {\bf B} \times {\bf A} \Big)||}\ ,
\end{equation}
provided  that  ${\bf  A}  \times  {\bf B}  \neq  {\bf  0}$,  i.e.~the
3-vectors ${\bf A}$ and ${\bf B}$  are not parallel to each other.  If
${\bf  A}$ and ${\bf  B}$ are  parallel, the  solution is  not unique,
given by ${\bf \widehat{x}} =  {\bf A} \times {\bf N}/||{\bf A} \times
{\bf N}||$, where  ${\bf N}$ is an arbitrary  3-vector non-parallel to
${\bf A}$  or ${\bf B}$, i.e.~${\bf  N}\cdot {\bf A} \neq  0$. In this
case,  the  unit-norm  vector  ${\bf  \widehat{x}}$ lies  on  a  plane
perpendicular to both ${\bf A}$  and ${\bf B}$, and ${\bf B}\cdot {\bf
\widehat{x}} = 0$.

If ${\bf A}$ and ${\bf B}$ are not parallel to each other and the norm
of the  vector ${\bf x}(t)$ is fixed  to be $||{\bf x}(t_*)||  = R$ at
$t=t_*  > 0$,  then the  largest value  of ${\bf  B}\cdot {\bf  x}$ at
radius $R$ is determined by
\begin{equation}
 \label{Ooptimal}
{\bf B}\cdot {\bf x}(t_*)\ =\ t_*\; ||{\bf A}\times {\bf B}||^2\ =\ 
t_*\; \Big[\, ||{\bf A}||^2\, ||{\bf B}||^2\: -\: 
({\bf A}\cdot {\bf B})^2\,\Big]\ .
\end{equation}
Likewise, the minimum  of ${\bf B}\cdot {\bf x}$  at radius $R$ occurs
at $t=-t_*$.

\section{A Higher-Dimensional Analytic Construction}

The above geometric construction for the 3D example can be generalized
to  higher dimensions,  where  ${\bf B}$  is  an $n$-dimensional  real
vector and $[{\bf A}]$ is an $m\times n$-dimensional real matrix (with
$0\leq  m<n$) that imposes  a number  $m$ of  real constraints  on the
$n$-dimensional  real vector  ${\bf x}$  through the  condition $[{\bf
A}]\cdot {\bf x}  = {\bf 0}$, with ${\bf 0}$ being  the null vector in
$m$ dimensions.

Our analytic construction proceeds along  the lines of the 3D example.
In  detail, in  analogy to  the  3D vector  ${\bf A}$,  we define  the
$m$-form  ${\bf  A}^{(m)}$  through  the $m$-fold  exterior  or  wedge
product:
\begin{equation}
  \label{Amform}
({\bf A}^{(m)})_{\alpha_1 \alpha_2 \dots \alpha_m}\ =\ 
({\bf A}_1 \wedge {\bf A}_2 \wedge
\dots  {\bf A}_m )_{\alpha_1 \alpha_2 \dots \alpha_m}\ =\ 
                A_{[1\alpha_1}\,A_{2\alpha_2}\, \dots\, A_{m\alpha_m ]}\; ,
\end{equation}
where the Greek indices  $\alpha_{1,2,\dots,m}\ =\ 1,2,\dots, n$ label
the columns of  the $m\times n$-dimensional matrix $[{\bf  A}] = ({\bf
  A}_1, {\bf A}_2, \dots, {\bf  A}_m )^T$.  The square brackets on the
RHS  of~(\ref{Amform})  indicate   that  the  $m$-rank  tensor  $({\bf
  A}^{(m)})_{\alpha_1  \alpha_2  \dots  \alpha_m}  \equiv  A_{\alpha_1
  \alpha_2  \dots \alpha_m}$  is  fully antisymmetrized  in all  Greek
indices, entailing that ${\bf A}^{(m)}$ is an $m$-form.

Correspondingly,  the analogue  of the  direction ${\bf  A}\times {\bf
  B}$, which determines  the normal to the plane  defined by ${\bf A}$
  and ${\bf B}$ in the 3D example of Section~2, is the $(n-m-1)$-form:
\begin{equation}
  \label{Bform}
({\bf C}^{(n-m-1)})_{\gamma_1 \gamma_2 \dots \gamma_{(n-m-1)}}\ =\ 
\varepsilon_{\alpha_1\alpha_2\dots \alpha_m \beta \gamma_1 \gamma_2
  \dots \gamma_{(n-m-1)}}\,  
A_{\alpha_1 \alpha_2 \dots \alpha_m}\, B_\beta\; ,
\end{equation}
where    summation   over    repeated   indices    is    implied   and
$\varepsilon_{\alpha_1\alpha_2\dots     \alpha_m     \beta    \gamma_1
\gamma_2\dots\gamma_{(n-m-1)}}$  is   the  usual  fully  antisymmetric
Levi--Civita tensor generalized to $n$ dimensions.  In the language of
differential  forms,  $({\bf  C}^{(n-m-1)})_{\gamma_1  \gamma_2  \dots
\gamma_{(n-m-1)}}\     \equiv\      C_{\gamma_1     \gamma_2     \dots
\gamma_{(n-m-1)}}$  is,  up  to  an  irrelevant  overall  factor,  the
Hodge-dual  product between the  1-form $({\bf  B}^{(1)})_\beta \equiv
B_\beta$,  representing  the components  of  the $n$-dimensional  real
vector ${\bf B}$, and the $m$-form ${\bf A}^{(m)}$, i.e.,
\begin{equation}
  \label{Cform}
{\bf C}^{(n-m-1)} \ =\ \star \Big({\bf B}^{(1)} \wedge {\bf A}^{(m)}\Big)\; ,
\end{equation}
where $\star$ denotes the standard Hodge-star operation~\cite{SG,SMC}
applied within the entire $n$-dimensional vector space.

The components  $x_\alpha$ of the optimal  $n$-dimensional real vector
${\bf x}$  are given by the  Hodge-dual product of  the $m$-form $({\bf
  A}^{(m)})_{\alpha_1  \alpha_2  \dots  \alpha_m}  \equiv  A_{\alpha_1
  \alpha_2   \dots   \alpha_m}$    and   the   $(n-m-1)$-form   $({\bf
  C}^{(n-m-1)})_{\gamma_1  \gamma_2   \dots  \gamma_{(n-m-1)}}  \equiv
C_{\gamma_1 \gamma_2 \dots \gamma_{(n-m-1)}}$, which is the 1-form:
\begin{eqnarray}
  \label{xnD}
x_\alpha(t) \! &= &\! t'\,
\varepsilon_{\alpha \alpha_1\alpha_2\dots \alpha_m \gamma_1 \gamma_2
  \dots \gamma_{(n-m-1)}}\, 
A_{\alpha_1 \alpha_2 \dots \alpha_m}\; C_{\gamma_1 \gamma_2 \dots
  \gamma_{(n-m-1)}}\nonumber\\ 
\! & =&\! t\,  
\varepsilon_{\alpha \alpha_1\alpha_2\dots \alpha_m \gamma_1 \gamma_2
  \dots \gamma_{(n-m-1)}}\,
\varepsilon_{\beta \beta_1\beta_2\dots \beta_m\gamma_1 \gamma_2 \dots
  \gamma_{(n-m-1)}}\nonumber\\ 
\! &&\! \times\,
A_{1\alpha_1}\,A_{2\alpha_2}\, \dots\, A_{m\alpha_m}\, B_\beta\, 
A_{1\beta_1}\,A_{2\beta_2}\, \dots\, A_{m\beta_m}\ ,
\end{eqnarray}
where $t$ and $t'$ are arbitrary real parameters. Equation~(\ref{xnD})
represents  the central  result of  this note, and can be equivalently 
cast into the more compact form:
\begin{equation}
  \label{xForm}
{\bf x} (t)\ =\ t'\; \star\!\Big({\bf A}^{(m)} \wedge {\bf
  C}^{(n-m-1)}\Big)
\ =\
 t\; \star\!\Big[ {\bf A}^{(m)} \wedge\, 
\star \Big({\bf B}^{(1)} \wedge {\bf A}^{(m)}\Big)\Big]\; .
\end{equation}
We emphasize  the complete analogy between the  general solution given
by the RHS of (\ref{xForm}) and the 3D result of (\ref{Param3}).

By construction, the  $n$-dimensional 1-form vector ${\bf x}_\alpha(t)
\equiv  x_\alpha  (t)$ satisfies  the  linear  system of  constraints:
$[{\bf A}]  \cdot {\bf  x} =  {\bf 0}$.  Within  a fixed  given radius
$||{\bf x}(t_*)||  = R$  at $t=t_*  > 0$, the  maximum value  of ${\bf
  B}\cdot {\bf x} (t_*)$ is given by
\begin{eqnarray}
  \label{Optimaln}
{\bf B}\cdot {\bf x} (t_*) \!& = &\! t_*\;
\varepsilon_{\alpha \alpha_1\alpha_2\dots \alpha_m \gamma_1 \gamma_2
  \dots \gamma_{(n-m-1)}}\,
\varepsilon_{\beta \beta_1\beta_2\dots \beta_m \gamma_1 \gamma_2 \dots
  \gamma_{(n-m-1)}}\nonumber\\ 
\! &&\! \times\, B_\alpha
A_{1\alpha_1}\,A_{2\alpha_2}\, \dots\, A_{m\alpha_m}\, B_\beta\, 
A_{1\beta_1}\,A_{2\beta_2}\, \dots\, A_{m\beta_m}\ .
\end{eqnarray}
The  minimum value is obtained correspondingly as $t =  -t_*$. Obviously, a
non-zero value for  ${\bf B}\cdot {\bf x} (t_*)$  is obtained, iff the
$n$-dimensional vectors $A_{1\alpha_1},  \, A_{2\alpha_2}, \, \dots,\,
A_{m\alpha_m}$  and $B_\beta$  are  all linearly  independent of  each
other.

\noindent
{\bf  Proof.}   We  now  present   a  simple  proof  of  the  analytic
construction given in~Eq.~(\ref{Optimaln}).   To this end, we consider
the    non-trivial   case,    where   the    $n$-dimensional   vectors
$A_{1\alpha_1}$, $A_{2\alpha_2},$ \dots, $A_{m\alpha_m}$ and $B_\beta$
span a  non-degenerate $(m+1)$-dimensional subspace.   Consequently, a
linear  transformation represented  by the  $m\times m$  matrix $[{\bf
R}]$  can be  performed on  the  left of  the $m\times  n$-dimensional
matrix $[{\bf A}]$, i.e.~$[{\bf A'}] = [{\bf R}]\cdot [{\bf A}]$, such
that  $[{\bf  A'}]\cdot {\bf  x}(t_*)  =  {\bf  0}$ and  $A'_{k\alpha}
A'_{l\alpha} = \mbox{const.}   \delta_{kl}$, with $k,l = 1,2,\dots,m$.
Specifically,  this   linear  transformation  exemplifies   the  usual
Gram--Schmidt  approach~\cite{SL}  to  orthogonalizing  a set  of  $m$
linearly-independent  vectors  which  span  the  same  $m$-dimensional
subspace  as   the  original  vectors.   In   addition,  the  complete
$n$-dimensional space can be  rotated by an orthogonal transformation,
i.e.,~${\bf x}(t_*) \to {\bf  x}'(t_*) = [{\bf O}]\cdot {\bf x}(t_*)$,
without affecting  the norm of all $n$-dimensional  vectors, such that
$A''_{k\alpha}  =  a_k \delta_{k\alpha}$.   Here,  $a_k$ are  positive
constants  that give  the norms  of the  orthogonalized  vectors ${\bf
A'}_k$ which are equal to those of the orthogonally-rotated ones ${\bf
A''}_k$.  In this linear basis, only the components of $B_\beta$ lying
in  the complementary  dimensions $m+1,m+2,\dots,  n$ give  a non-zero
contribution in Eq.~(\ref{Optimaln}).   If we denote this orthogonally
projected  vector with  $B^\perp_\beta$, the  linear  expression ${\bf
B}\cdot {\bf x} (t_*)$ becomes
\begin{equation}
{\bf B}\cdot {\bf x} (t_*) \ =\ t_* \prod_{k=1}^m a^2_k\ \sum_{\alpha=
  m+1}^n (B^\perp_\alpha)^2\ .
\end{equation}
This  represents  indeed  the  largest  value  for  the  above  linear
expression.  Note that the  optimal direction ${\bf x}(t)$ is parallel
to the reduced vector: ${\bf B}^\perp$.

\section{Conclusions}

The  problem studied  in  this note  deals  with a  specific class  of
problems within the wider context  of linear programming.  It has been
shown that the optimal  unit-norm direction ${\bf \widehat{x}}$ of the
objective  ${\bf B}\cdot {\bf  \widehat{x}}$, subject  to a  system of
linear constraints of the form $[{\bf  A}]\cdot {\bf x} = {\bf 0}$, is
given by~(\ref{xnD}), where~$[{\bf  A}]$ is an $m\times n$-dimensional
real matrix.  In particular, the  optimal solution can be expressed in
terms  of wedge  and  Hodge-dual products  of  differential forms,  as
stated  in~(\ref{xForm}).  The latter  provides an  explicit geometric
connection  of  the  3D  result  in~(\ref{Param3})  with  the  general
$n$-dimensional solution given in~(\ref{xForm}).

A straightforward extension of the problem studied in the present note
would  be   to  consider   a  complexification  of   the  constraints,
i.e.,~$[{\bf  A}]$ is  an $m\times  n$ complex  matrix, in  which case
${\bf  x}$ is  a complex  $n$-dimensional vector.   Such  an extension
makes  sense if  the objective  is a  real number,  e.g.,~the  real or
imaginary  part of  the expression  ${\bf B}\cdot  {\bf \widehat{x}}$.
This problem  can be solved along  the lines of  the presented method,
where  ${\bf x}$  is  treated  as $2  n$-dimensional  real vector  and
turning the  complex $m\times n$  matrix $[{\bf A}]$ into  a $2m\times
2n$ real matrix.

The present  geometric solution  may suggest the  existence of  a more
profound connection between differential forms and linear programming,
beyond the  limitations of the  specific problem under study.   It may
therefore  be interesting  to  pursue this  direction of  mathematical
research  in  linear  programming   with  greater  vigour  in  future.
Likewise,   it  would   be  interesting   to  explore   the  potential
applications  of this  geometric solution  to a  plethora  of problems
related  to  physics~\cite{ELPouter},  computer science,  engineering,
biology, actuarial science and economics.

\subsection*{Acknowledgements}

We thank  Luis Alvarez-Gaum{\'e} and Nikolaos  Papadopoulos for useful
comments and  suggestions, as  well as Frank  Deppisch for  a critical
reading of the note.

\end{document}